\documentclass[11pt]{article}
\usepackage[utf8]{inputenc}
\usepackage{url,geometry}
\usepackage{bbm}
\usepackage{bm}
\usepackage{amsmath, scalerel, amsthm, amssymb}
\usepackage{lscape}
\usepackage{algorithm}
\usepackage{algorithmic}
\usepackage{graphics}
\usepackage{graphicx}
\usepackage{wrapfig}
\usepackage{multirow,bigdelim}
\usepackage{pict2e,color}
\usepackage{amsfonts}
\usepackage{verbatim}
\usepackage{booktabs}
\usepackage{url}
\usepackage{longtable}
\usepackage{tabularx}
\usepackage{multicol}
\usepackage{enumerate}
\usepackage{stmaryrd}
\usepackage[normalem]{ulem}
\usepackage{caption}
\usepackage{subcaption}
\usepackage{float}
\usepackage{tikz}
\usepackage{pgfplots}
\usepackage{setspace}
\usepackage{xurl}
\newcommand{\exclude}[1]{}
\usepackage{nicefrac}
\usepackage{adjustbox}
\usepackage[round]{natbib}
\usepackage{todonotes}

\newtheorem{thm}{Theorem}

\DeclareMathOperator*{\argmax}{arg\,max}

\newcommand{\added}[1]{\textcolor{black}{#1}}

\title{\added{Interdicting attack plans with boundedly-rational players and multiple attackers}: An adversarial risk analysis approach}

\author{\textbf{Eric DuBois, Ashley Peper, and Laura A.~Albert\footnote{Corresponding author}}\\
University of Wisconsin-Madison\\
Industrial and Systems Engineering \\
1513 University Avenue \\
Madison, Wisconsin 53706 \\
Email: \url{laura@engr.wisc.edu} 
}

\date{February 2023}

\begin{document}
\maketitle

\begin{abstract}
\added{Cybersecurity planning} supports the selection of and implementation of security controls in resource-constrained settings to manage risk. 
\added{Doing so requires considering} adaptive adversaries with different levels of strategic sophistication in modeling efforts to support risk management.  However, most models in the literature only consider rational \added{or non-strategic} adversaries.
Therefore, we study how to inform defensive decision-making to mitigate the risk from \added{boundedly rational players, with a particular focus on making integrated, interdependent planning decisions}.
\added{To achieve this goal, we introduce a modeling framework for selecting a portfolio of security mitigations that interdict adversarial attack plans that uses a structured approach for risk analysis.}
\added{Our approach adapts} adversarial risk analysis and cognitive hierarchy theory to consider a maximum reliability path interdiction problem with a single defender and multiple attackers who have different goals and levels of strategic sophistication. 
\added{Instead of enumerating all possible attacks and defenses, we} introduce a solution technique based on integer programming and approximation algorithms to iteratively solve the defender's and attackers' problems. 
A case study illustrates the proposed models and \added{provides insights into defensive planning}.
\end{abstract}


\section{Introduction}
Cybersecurity is an important concern for governments and organizations throughout the world due to the growing reliance on digital connectivity and the growing number of threats. Cyber attacks are increasingly common, costing the U.S. economy \$57-\$109 billion in 2016 \citep{council_of_economic_advisors_cost_2018} and affecting systems throughout the economy, including in healthcare \citep{kruse_cybersecurity_2017}, energy \citep{wang2013cyber}, and industrial control \citep{knowles_survey_2015}. Many possible security controls exist to mitigate these risks \citep{NIST800}.
\added{Cybersecurity planning requires periodically selecting a portfolio of security controls (e.g., on an annual basis), which allows an organization to manage the risk associated with  vulnerabilities that have emerged.}
\added{However, many organizations find it challenging to keep up with selecting and deploying security controls given that they operate in resource-constrained environments \citep{stevens2020compliance}. }

\added{There is a growing body of literature that applies risk analysis techniques to manage cybersecurity risk through the strategic prioritization of security controls. 
Some of these efforts model attackers as non-strategic players using probability distributions and prioritize security controls in rank order based on their cost-effectiveness \citep{hubbard2016measure}. 
Increasingly, security controls are prioritized using a structured approach to aid in planning efforts with well-defined goals and threat scenarios (NIST \citeyear{nist2018}), where attack graphs are used to represent known vulnerabilities and visualize potential mitigations \citep{lallie2020review}.
%
Recent research uses integer programming \citep{zheng_budgeted_2019} and robust optimization \citep{zheng2019robust} to select security mitigations using a structured approach based on attack graphs, however, these papers assume attackers are either not strategic or limited to selecting a worst-case scenario. 
}


\added{A stream of papers in the literature explicitly consider adaptive adversaries in cybersecurity planning through the application of adversarial risk analysis (ARA). In these models, a single defender, the security planner, selects security controls that perform well given that adaptive adversaries will attempt to work around any new security controls that are put in place. 
ARA frameworks are versatile and have motivated defender-attacker models that capture a wide range of conditions and assumptions \citep{banks2020adversarial}.
%
As a result, ARA has been applied to many of these defender-attacker models in various security settings \citep{rios2012adversarial} and has been adapted to cybersecurity models with multiple adversaries with different levels of intentionality \citep{insua_adversarial_2019}. 
A limitation of ARA approaches is that they enumerate all possible attacks and defenses \citep{banks2020adversarial}, with \cite{wang2011network} as an exception. Thus, ARA algorithms are intractable when it is not practical to enumerate cybersecurity attack and defense choices. 
\citet{zheng2019interdiction} seek to overcome this challenge by introducing a structure that simplifies computational requirements based on a network interdiction model that delays the attack plans of multiple attackers. The structure introduced by the network interdiction model allows for the use of integer programming algorithms to solve for the defender and attacker strategies. 
}

%
\added{We build upon previous work by introducing an ARA framework that considers boundedly rational players, including a defender and multiple adversaries, to inform the selection of security controls that interdict adversarial attack plans. 
While our approach is motivated by cybersecurity planning, it can be used more broadly in the security context where players seek to maximize or minimize the probability of attack in a multi-layer defense system.
}


\subsection{Approach}\label{subsec:approach}
In the cybersecurity planning problem we consider, a defender seeks to minimize the probability of a successful attack by multiple attackers with different levels of strategic sophistication over a planning horizon by selecting a portfolio of security controls subject to a budget.
%
%
%
We adapt an ARA framework \citep{insua_adversarial_2019} to capture the strategic selection of a portfolio of security controls given that the defender and the attackers are boundedly rational. 
This modeling approach allows us to inform defensive decisions and planning against a range of attackers, which more accurately reflects the system we are modeling \citep{scheibehenne_can_2010}.
\added{Since new vulnerabilities emerge on a regular basis, cybersecurity planning should be performed regularly (e.g., annually), and the methods in this paper can aid in this process.}

\added{Attack modeling is an important step in cybersecurity planning. In vulnerability analysis, vulnerabilities can be characterized by various steps required to successfully carry out an attack \citep{Schneier1999}, which provides a structured approach to represent attack scenarios. In a graph structure, the nodes represent attack states (e.g., the choice of attack-type, the target of the attack, or attack milestones), edges represent intermediate exploits in an attack. The difficulty of an adversary completing an exploit is captured by a conditional probability of successfully traversing an arc. 
A path from root to leaf corresponds to an attack against the system, and therefore, we view an attack as a path in a graph between a source and sink node.   Given that there are many adversaries who have different knowledge of the vulnerabilities and different capabilities, the attack graphs may have topology and parameters specific to each adversary.}

\added{An adversary's probability of successfully carrying out all exploits in an attack is captured by the probability they traverse the network on the path they select. Security controls interact with attack graphs by decreasing the probability that the completion of individual exploits (arcs that are traversed). Security controls may encourage adversaries to select alternative paths.
The defender uses their private information and the paths in the attack graphs which they believe the attacker will choose to determine their choice of security controls. }

\added{These strategic interactions motivate the application of the maximum reliability path interdiction problem to the planning problem under consideration. Network interdiction models, including the maximum reliability path interdiction, have been widely used to model attacker-defender games, usually assuming two players and rational decision-makers \citep{smith_survey_2019}. Although some papers have lifted the assumptions of shared information (e.g., \cite{salmeron_deception_2012}) and shared beliefs regarding the probabilities of traversing edges in the network (e.g., \cite{morton_models_2007}), none have considered boundedly rational players. This paper seeks to fill this gap.
}



We adapt and apply ideas from the adversarial risk analysis (ARA) framework presented by \citet{insua_adversarial_2019} to the maximum-reliability network interdiction problem to support decision-making for cybersecurity planning.
%
An ARA approach allows us to consider adversaries that are not rational and who may have varying levels of strategic sophistication, important features of the application under consideration. In particular, we model players who are boundedly rational using cognitive hierarchy theory \citep{camerer_cognitive_2004} and level $k$ thinking \citep{stahl_players_1995}. 

\begin{wrapfigure}{r}{0.5\textwidth}
    \centering
    \includegraphics[width=7cm]{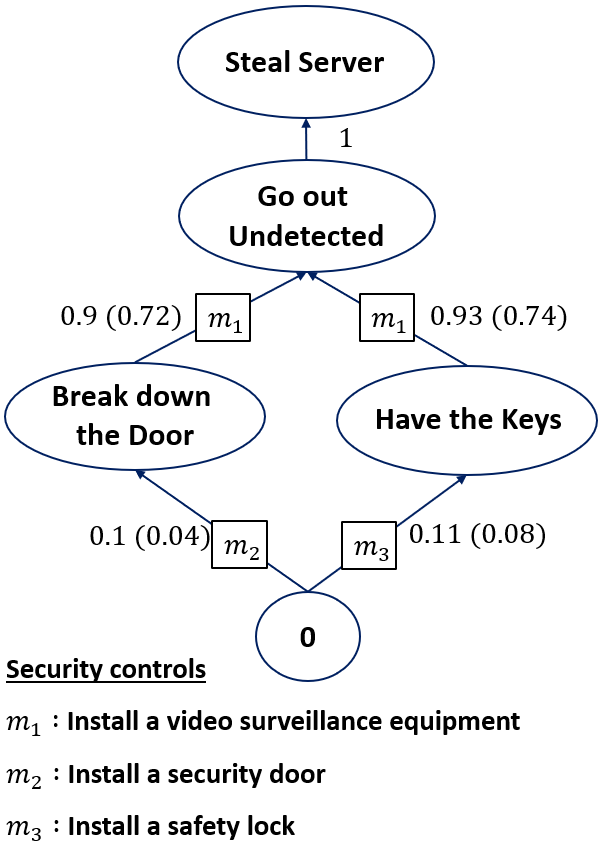}
    \caption{An illustrative example of a single attack graph with three security controls}
    \label{fig:example}
\end{wrapfigure}

\added{To illustrate the approach taken in this paper, consider the following example with a single defender and two types of adversaries. In the example, adapted from \cite{bistarelli2006defense}, adversaries  attempt to steal a server. Figure \ref{fig:example} captures the attack graph under consideration, including the exploits, and two pathways to steal the server (reach the top node). There are three security controls ($m_1, m_2$, and $m_3$) that are listed along the arcs they interdict. For illustration purposes we assume the defender can select one control. The un-interdicted traversal probabilities are listed next to each arc, and the interdicted traversal probability is listed in parentheses when the security control is in use. Without any security controls, the left and right pathways to the ``Steal server'' node have traversal probabilities of 0.09 and 0.102, respectively. }

\added{In this example, an opportunistic attacker with a low degree of strategic sophistication who ignores possible security defenses selects the right path, since it has a higher un-interdicted traversal probability. A defender who anticipates only attacks from this type of attacker would select $m_3$, since this mitigation lowers the traversal probability of the right path the most. A slightly more strategic attacker would anticipate $m_3$ being deployed and would then choose the left path. This would lead a more strategic defender to select control $m_2$ to defend against attackers attempting only the left path. A defender who anticipates both types of attackers, 25\% of whom are somewhat strategic and select the left path and 75\% of whom are opportunistic and attempt the right path, would select mitigation $m_1$. 
}

\added{This simple example based on a single attack graph highlights how a methodology that considers multiple, boundedly rational attackers can be informative for defensive planning decisions. In general, there are many attack graphs that capture various attack vectors as well as adversarial goals and capabilities. The security controls can be specific and delay a single exploit, such as $m_2$ in the previous example, or many exploits if the controls are cross-cutting, such as deploying multi-factor authentication or an employee training program. \cite{singhal2017security} and \cite{lallie2020review} provide additional guidance surrounding how to model attack graphs. }

\subsection{Contributions}
In summary, this paper makes the following contributions:
\begin{itemize}
\item We formulate the security control investment problem in an ARA framework as a maximum-reliability path interdiction game between a single defender and multiple attackers, \added{all of whom are boundedly rational} with differing levels of strategic sophistication.  The defender and attackers' problems are formulated as mixed integer programming models. 
\item We introduce solution techniques based on mixed integer programming algorithms and approximation algorithms. The defender and attackers' problems are solved iteratively, and the inputs to each model are updated after each iteration. We prove that the defender's problem is equivalent to a submodular maximization problem subject to a budget constraint, which enables the use of heuristics for identifying solutions that are at least $1-1 / e$ of the optimal solution value.
\item We apply the modeling approach to a case study. We identify solutions across a range of levels of strategic sophistication and consider the effects of the defender misjudging the sophistication of the attackers. 
\end{itemize}

The organization of this paper is as follows. We first survey the literature in Section \ref{cyber:sec:litreview}.  In Section \ref{cyber:sec:ARA} we describe the ARA framework. Section \ref{cyber:sec:model} introduces the mixed integer programming models that capture the maximum-reliability path interdiction from the attacker and the defender perspectives. We introduce the approximation algorithms in Section \ref{sec:algorithms}, and we describe the ARA algorithm in Section \ref{sec:araFramework}. We present and analyze the case study and other computational results in Section \ref{cyber:sec:computational}. Section \ref{cyber:sec:conclusion} contains concluding remarks.

\section{Literature Review}
\label{cyber:sec:litreview}
The topics in this paper are related to several different areas of research, including adversarial risk analysis, boundedly rational thinking, security control investment, and maximum-reliability path interdiction. We present a summary of the most relevant papers in the literature. 

\subsection{Adversarial Risk Analysis}
\added{ARA has been widely applied to security applications \citep{rios2012adversarial}. 
Several papers use ARA framework for cybersecurity, focusing on defender-attacker games \citep{wang2019adversarial}, insider threat modeling \citep{joshi2020insider}, and adversarial machine learning and data manipulation \citep{caballero2021adversarial}. 
\citet{insua_adversarial_2019} apply the ARA framework to a cybersecurity resource allocation problem to inform a portfolio of defensive actions, including the purchase of cyber-insurance, by including both intentional and non-intentional threats. 
%
%
} 

As with most ARA models, \citet{insua_adversarial_2019} enumerate all possible attacks and defenses. \citet{cano_modeling_2016} similarly apply ARA to a cybersecurity setting by enumerating specific attacks and defenses to determine an optimal security allocation to minimize disruptions to airport operations. To our knowledge, only \citet{wang2011network} consider an ARA framework in which the attacks and defenses are not enumerated. They consider the optimal path for a convoy through a network where an attacker has placed several improvised explosive devices (IEDs) at nodes within the network. The defender seeks to minimize the routing cost. By using a network model, they compactly represent many possible convoy routes. Since \citet{wang2011network} utilize the additive nature of their cost function to efficiently solve their problem, it is not possible to apply their solution method to our problem. \added{In contrast to the existing literature, we consider interdicting attack plans to support cybersecurity planning, and we introduce a methodology to solve defender and attacker problems based on integer programming and approximation algorithms, since enumerating the attack and defense choices is intractable.}



\subsection{Bounded Rationality}
We build on the work of \citet{rothschild2012adversarial}, who develop an algorithm for applying bounded rationality, specifically level $k$ thinking, within an ARA framework. Level $k$ thinking begins with non-strategic level 0 thinkers who act without regard to other players. 
In comparison to level $k$ thinking, where a level $k$ player optimizes over only the level $k-1$ opponent, cognitive hierarchy theory assumes that a level $k$ thinker optimizes over a distribution of players between level 0 and level $k-1$ \citep{camerer_cognitive_2004}. We use this method to model the defender due to the multiple-attacker scenarios they face. While the logic of level $k$ thinking and cognitive hierarchy theory are theoretically subject to infinite regression, empirically it has been found that most people do not think beyond level 2 or level 3 \citep{lee_game_2012}.
Higher-level thinkers ($k>0$) assume that their opponent is a level $k-1$ thinker, e.g., level 1 players optimize against level 0 players, level 2 players optimize against level 1 players, and so on. \citet{rothschild2012adversarial} create an algorithm for determining the strategy of a level $k$ opponent by using recursion to build belief distributions regarding the attacks or defenses that the player uses. \added{Considering bounded rationality in the adversaries is a novel aspect of our paper. To the best of our knowledge, our paper is the first to consider boundedly rational players when considering interdependent defensive decisions as considered in network interdiction. } 


\subsection{Security Control Portfolio Selection}
\added{A stream of the literature studies how to invest in  security controls given a limited budget \citep{fielder_decision_2016,zheng2019robust}. 
Many models consider non-strategic attacks and do not necessarily select controls that adequately protect against adaptive adversaries \citep{zheng_budgeted_2019}. Non-strategic attacks (also called opportunistic or non-targeted) continue to be carried out in roughly the same manner and with the frequency regardless of the defender’s security decisions. For example, models using a decision theory framework assume that the defender’s decision has no impact on the frequency of each method of attack \citep{cavusoglu_decision-theoretic_2008}. While some research seeks to allocate a budget in the presence of strategic attackers and natural disasters \citep{zhuang2007balancing}, none of the papers in this area consider how to select a portfolio of security defenses with multiple, boundedly rational adversaries.} 

\subsection{Maximum-Reliability Path Interdiction}
\added{Network interdiction models have been widely applied to infrastructure protection and resilience problems, where they inform how to protect critical components in a system to reduce worst-case vulnerability. Network interdiction problems are modeled as Stackelberg or Cournot game models of defender-attacker games,} usually assuming two players and rational decision-makers.
\citet{smith_survey_2019} provide a recent survey of this area. A shortest path interdiction problem may be used to model attacks or projects where the attacker seeks to minimize the time required to complete an attack (traverse the network) and the defender seeks to maximally delay an attack by interdicting (lengthening) edges on the network \citep{israeli_shortest-path_2002}. The shortest path interdiction problem is an equivalent formulation to the maximum-reliability interdiction problem \citep{morton_models_2007}.
Network interdiction models have been extended to include imperfect and private information \citep{salmeron_deception_2012}. \added{Several researchers have recommended that network interdiction models be extended to consider boundedly rational players \citep{zhang2018stochastic} and other more realistic features to aid in defensive planning \citep{albert2023homeland}.
}


\added{A stream of papers study how to interdict attack graphs to inform defensive cybersecurity planning efforts. \cite{nandi2016interdicting} introduce a bi-level defender-attacker model to help organizations select and deploy security countermeasures by interdicting attack graphs.  \citet{Letchford-Vorobeychik2013} introduce a different Stackelberg game in which a defender seeks to interdict an attack plan, lifting the assumption that the game is zero-sum. 
\citet{zheng2019interdiction} introduce a bi-level network interdiction model that seeks to identify a portfolio of security controls that maximially delay a large number of adversarial cyber attacks from multiple attackers under uncertainty.
}

\added{In sum, our paper adds to the literature by combining ARA and cognitive hierarchy theory with network interdiction modeling to inform defensive cybersecurity planning.}



\section{Adversarial Risk Analysis Framework}
\label{cyber:sec:ARA}
In this section, we introduce the ARA framework in this paper as well as how we model the boundedly rational players.  
ARA takes the perspective of one player, the defender in this paper, and seeks the optimal action for that player based on the actions/reactions they believe the other players will take.

\subsection{Player Interactions in the Game}
The maximum reliability network interdiction modeling approach considers a one-off encounter between a single defender and multiple attackers.  This is reasonable for cybersecurity investment decisions considered in the case study, where planners must protect information systems from many attacks. In our approach, attackers do not know which security controls have been selected by the defender. The attackers seek to maximize the probability of their attacks succeeding based on perceived defenses, and the defender seeks to prevent an attack. Beliefs about the defender's budget and security control costs are finite and discrete. 

We now describe the general form of the game that we consider.  
There are a set of attackers who each choose a method of attack and target that maximizes their probability of success, which defines their \textit{reliability}. This is equivalent to choosing a path through the network. By choosing a path, the attacker chooses the set of edges they traverse that determines the path's reliability. By selecting controls, the defender reduces the reliability for each of the edges that the controls cover. The defender selects a portfolio of controls to maximize the conditional probability that an attack is prevented subject to a budget constraint.

Most of the information regarding the structure of this network is treated as shared beliefs, whereas the parameters of the network, such as the edge reliabilities, are private beliefs.  Specifically, we assume that the set of nodes, edges, and possible controls are shared between all players. This is equivalent to assuming that the defender knows all possible threats the attackers may consider. \added{The approach informs defensive planning based on known vulnerabilities, and it can be updated to include new vulnerabilities that have been discovered. }
Other information in the game may be private.  This includes the edge reliabilities, the effectiveness of controls at decreasing those reliabilities, the proportion of attacks from each attacker, the cost to institute a control, and the defender's budget.

\subsection{\label{cyber:sec:bounded_rational}Boundedly Rational Players}
From the defender's perspective, the attackers may have different levels of strategic sophistication. Therefore, a level $k$ defender plans for level $0, 1, ...,k-1$ attackers. An attacker's strategy represents a path.  The defender's strategy is a portfolio that consists of a set of controls whose total cost is within the budget.

Since we are seeking an optimal defense, we find strategies for level $k$ defenders up to maximum level $K$---the defender's level---and strategies for level $k$ attackers up to maximum level $K-1$. Lower level attacker and defender strategies are recorded as the algorithm progresses, which provides a suite of portfolio options corresponding to the differing defender levels of strategic sophistication. By analyzing the various portfolios within this suite, we can make more informed decisions about how the defender's posture should change when confronted by more or less ``sophisticated'' attackers.

The basic algorithm to compute the strategies of different level $k$ attackers of each type is as follows. We elaborate upon this algorithm in detail in Section \ref{sec:algorithms}. We begin with information and beliefs about the system. 
We start with $k=0$ and construct level $0$ attacker and defender strategies for nonstrategic players. We then increment the value of $k$ by $1$ and compute a level $k$ defender's optimal portfolio using the defender's program, \textit{OptDef}, presented in Section \ref{cyber:sec:optdef}. If the strategy for the highest level of defender, $K$, has been calculated, the algorithm terminates. Otherwise, we calculate the level $k$ attacker's optimal attack path using the attacker's program, \textit{OptAtt}, presented in Section \ref{cyber:sec:optatt}. We repeat this step by incrementing $k$ by $1$ and solving the defender's and attackers' problems until the algorithm terminates.

\section{\label{cyber:sec:model}Model Formulations}

In this section, we formulate a simultaneous single-defender multiple-attacker game based on a maximum reliability network interdiction problem.  We do so by introducing optimization problems from the attackers' perspective, \textit{OptAtt}, and the defender's perspective, \textit{OptDef}. 
There is a single defender and a set of attackers $\mathcal{A}$.
Without loss of generality, each attacker begins at the super source node, $1$, and progresses through the graph to the super sink node, $n$. Each attacker chooses a path of maximum reliability, the probability that they believe they will successfully traverse the graph based on perceived defender decisions. 

Notation that reflects the defender's and attackers' beliefs have a $D$ and $A$ subscript/superscript, respectively. Beliefs that are updated based on the level of strategic sophistication have a $k$ superscript. Notation that captures decision variables does not explicitly include $k$, $D$, or $A$. 
Table \ref{tab:systemnotation} provides a summary of the notation relevant to the parameters and players beliefs as well as the decision variables. For clarity, we use Latin characters for variables and shared information, and we use Greek characters for beliefs. 

\begin{table}
\caption[Notation]{\label{tab:systemnotation}Notation notation}
\singlespacing
\begin{center}
\textbf{Sets: Common Information to Defender and Attackers }\\
\begin{tabular*}{\textwidth}{r c l}
\midrule
$\mathcal{A}$ & = & Attackers\\
$N$ & = & Nodes\\
$E$ & = & Edges\\
$E^+_i\left(E^-_i\right)\subseteq E$ & = & Edges that leave (enter) node $i\in N$\\
$M$ & = & Set of controls\\
$M_{ij}\subseteq M$ & = & Subset of controls that interdict edge $(i,j)\in E$\\
$P(i)$ & = & The set of paths from the source node to node $i\in N$ \\
\midrule
\\
\end{tabular*}


\textbf{Level $k$ Attacker $A \in \mathcal{A}$ Decisions}\\
\begin{tabular*}{\textwidth}{rcl}
\midrule
\multicolumn{3}{l}{\textbf{Attacker Variables}}\\
\multirow{2}{*}{$u_{ij}$} & \ldelim\{{2}{*}[=] &1, if edge $(i,j)\in E$ is on the attacker's path\\
&& 0, otherwise \\

$q_{ij}$ & = & probability that the attack of attacker $A \in \mathcal{A}$ reaches edge $(i,j)\in E$\\
\\
\multicolumn{3}{l}{\pmb{$\implies$}\textbf{ Level $k'$ (with $k'>k$) Defender Beliefs}}\\
$\psi_D^{k'}(i,j,A)$ & = & Conditional probability that a level $k' > k$ defender believes a level $k$ attacker \\ 
& &  $A \in \mathcal{A}$ attempts to traverse edge $(i,j)\in E^+_i$ given that they reach node $i$. \\
\midrule
\\
\end{tabular*}

\textbf{Defender $D$ Beliefs}\\
\begin{tabular*}{\textwidth}{r c l}
\midrule
$\delta_D(i,j,A)$ & = & Reliability of un-interdicted edge $(i,j) \in E$ for attacker $A \in \mathcal{A}$\\
$\tilde{\delta}_D(i,j,A)$ & = & Reliability of interdicted edge $(i,j) \in E$ for attacker $A \in \mathcal{A}$\\
$\theta_D(A)$ & = & Conditional probability that attacker $A$ attempts an attack \\
$\beta_D$ & = & defender's budget\\
$\kappa_D(m)$ & = &cost of control $m\in M$\\
\midrule
\\
\end{tabular*}

\textbf{Level $k$ Defender Decisions}\\
\begin{tabular*}{\textwidth}{rcl}
\midrule
\multicolumn{3}{l}{\textbf{Defender Variables}}\\
\multirow{2}{*}{$w_{m}$} & \ldelim\{{2}{*}[=] &1, if control $m\in M$ is chosen\\
&& 0, otherwise \\
\multirow{2}{*}{$x_{ij}$} & \ldelim\{{2}{*}[=] &1, if edge $(i,j)\in E$ is covered\\
&& 0, otherwise \\
$y_{ij}(A)$ & = & probability that attacker $A\in \mathcal{A}$ reaches un-interdicted edge $(i,j)\in E$\\
$\tilde{y}_{ij}(A)$ & = & probability attacker $A\in \mathcal{A}$ reaches interdicted edge $(i,j)\in E$\\
\\
\multicolumn{3}{l}{\pmb{$\implies$}\textbf{Level $k+1$ Attacker Beliefs}}\\
$\Omega_A^{k} = \{\omega_1,...,\omega_m\}$ & = & portfolio selected by a level $k$ defender \\
$\delta_A^{k+1}(i,j|\Omega^{k}_A)$ & = & Reliability of edge $(i,j) \in E$ under portfolio $\Omega_A^{k}$ for a level $k+1$ attacker  \\
\midrule
\\
\end{tabular*}

\end{center}
\end{table}

\subsection{\label{cyber:sec:optatt}Attackers' Problem}
We now consider the attackers' problem. 
We consider a directed acyclic graph $G=(N,E)$ consisting of a finite set of nodes $N$ and directed edges $E$.  The sets of edges leaving and entering a node $i \in N$ are denoted by $E^+_i$ and $E^-_i$, respectively.  Without loss of generality, we assume that the set of nodes $\{1,...,n\}\in N$ are ordered such that $i<j, \forall (i,j) \in E.$ 

A level 0 attacker is not strategic. We assume that a level 0 attacker $A$ employs a myopic greedy tactic for selecting a path. In this approach, edge reliabilities are given by the non-interdicted values, $\delta(i,j,A)$ for all $(i,j) \in E$, and each level 0 attacker selects the edge with the highest non-interdicted reliability from each node, breaking ties randomly. 
Other methods could be used to set level $0$ attacker paths.

We consider the reliability of an attack path for a level $k$ attacker.  A level $k$ attacker $A$ believes the level $k$-1 defender has chosen the portfolio $\Omega_A^{k-1} = \{\omega_1,...,\omega_m\}$. 
The attacker believes that the reliability of edge $(i,j)\in E$ is $\delta^k_A(i,j| \Omega_A^{k-1})$, its reliability with portfolio $\Omega_A^{k-1}$. We assume that the values of $\delta^k_A(i,j| \Omega_A^{k-1})$ are independent between edges.
Given $\Omega_A^{k-1}$, the attacker then believes the reliability of a fixed path $p$ to be $\prod_{(i,j)\in p}\delta^k_A(i,j|\Omega_A^{k-1})$. Each strategic attacker seeks to maximize the conditional probability that their attack succeeds. 
Using standard approaches based on recursion \citep{ahuja_network_1993}, we can represent the probability of an attack reaching edge $(i,j)$, captured by $q_{j \ell}$, using the pair of linear inequalities:
\begin{align*}
    q_{j \ell}\leq& u_{j \ell}\\
    q_{j \ell}\leq& \sum_{(i,j) \in E^-_j}\delta^k_A(i,j|\Omega_A^{k-1})q_{ij},
\end{align*}  
where the characteristic vector of the path chosen by the attacker is $u \in \{0,1\}^{|E|}$, where $u_{ij} = 1$ if edge $(i,j)$ is in the path, and $u_{ij} = 0$ otherwise. 


Using these expressions, we present the \textit{OptAtt} formulation for level $k$ attacker $A\in \mathcal{A}$ as a maximum reliability path problem:
\begin{align}
z_A^k = \max \; & \sum_{(i,n) \in E^-_n}\delta^k_A(i,j|\Omega_A^{k-1})q_{in}& \label{cyber:eqn:A0}\\
\label{cyber:eqn:A1} \text{s.t.}\; & \sum_{(1,i) \in E^+_1}u_{1i} \leq 1\\
& \label{cyber:eqn:A2} \sum_{(i,j) \in E^+_i}u_{ij} \leq \sum_{(j,i) \in E^-_i}u_{ji}  , &\forall i \in N\setminus\{1\}\\ 
& \label{cyber:eqn:A3} q_{ij} \leq u_{ij}, &\forall (i,j) \in E \\
& \label{cyber:eqn:A4} q_{j \ell} \leq \sum_{(i,j) \in E^-_j}\delta^k_A(i,j|\Omega_A^{k-1})q_{ij} , &\forall j \in N\setminus \{1\}, (j,\ell) \in E^+_j\\
& \label{cyber:eqn:A6} u_{ij} \in \{0,1\}, &\forall (i,j)\in E\\
& \label{cyber:eqn:A7} q_{ij} \geq 0, &\forall i \in N, (i,j)\in E
\end{align}

The objective \eqref{cyber:eqn:A0} for an attacker is to maximize the probability that their attack succeeds.  Constraints (\ref{cyber:eqn:A1}) -- (\ref{cyber:eqn:A2}) enforce that the $u_{ij}$ variables properly define a path.
Constraint (\ref{cyber:eqn:A1}) only allows the attacker to choose one attack path, and constraint set (\ref{cyber:eqn:A2}) preserves the balance of flow in and out of each node.  Constraint sets (\ref{cyber:eqn:A3}) -- (\ref{cyber:eqn:A4}) determine the probability that the attack succeeds.  Constraint set (\ref{cyber:eqn:A3}) allows attacks to progress only along the attack path.  Constraint set (\ref{cyber:eqn:A4}) balances the flow of attack probability in and out of each node. Constraint sets (\ref{cyber:eqn:A6}) and (\ref{cyber:eqn:A7}) require variables to take on binary and non-negative values, respectively.
\textit{OptAtt} is a canonical form maximum reliability path problem, which can be solved for each attacker as a shortest path problem with Dijkstra's algorithm after a negative logarithm transform of the edge reliabilities \citep{morton_models_2007}.

After solving \textit{OptAtt} for each attacker, we use the solution to construct the defender's beliefs about the paths that each attacker $A\in \mathcal{A}$ will take. Recall that a defender of level $k'$ defends against a set of attackers $\mathcal{A}$ that may have different levels of sophistication $k$, with $0 \le k < k'$. A solution to \textit{OptAtt} for a level $k$ attacker $A \in \mathcal{A}$ can therefore be used to inform the beliefs of the level $k'$ defenders, with $k+1 \le  k' \le K$. To do so, we derive the paths that a level $k$ attacker $A$ takes from the \textit{OptAtt} solution variables $u_{ij}$, which are converted to inform a level $k'$ defender's belief parameters. 
A level $k'$ defender believes the conditional probability attacker $A$ chooses edge $(i,j) \in E$ after reaching node $i$ is 
\begin{align}\label{eq:defenderupdate}
    \psi^{k'}_D(i,j,A) = u_{ij}, \ \forall (i,j) \in E .
\end{align}
While we do not explicitly model the level $k$ in the attacker variables in \textit{OptAtt}, the level of the attacker is implicitly retained when setting the values of $\psi^{k'}_D(i,j,A)$.


\subsection{\label{cyber:sec:optdef}Defender's Problem}
We now consider the defender's problem. 
%
The defender maximizes the probability that an attack is prevented by selecting controls from set $M$. Each control $m \in M$ has cost $\kappa_D(m)$ subject to defender budget $\beta_D$. Control $m$, interdicts a set of edges, decreasing the reliability of all edges in the set. The subset $M_{ij}\subseteq M$ contains all controls that interdict edge $(i,j)$.  This approach is slightly different than that of the canonical maximum reliability interdiction problem, where the defender directly interdicts edges on the network. We treat the budget and costs as beliefs, since the attacker may not know their values. We assume a level $0$ defender does not choose any controls.

The defender believes that the conditional probability that attacker $A \in \mathcal{A}$ attempts an attack is $\theta_D(A)$. This is assumed to be determined exogenously, possibly from expert opinion or a risk assessment.
The defender faces multiple attackers with different levels $k'$ with $0 \le k' \le k-1$, which are captured in the values of $\psi^{k}_D(i,j,A)$.
%
Each attack follows a fixed path $p$ with a given probability that depends on the defender's level. The defender has belief probabilities $\tilde{\delta}_D(i,j,A)$ and $\delta_D(i,j,A)$ that reflect the reliability of edge $(i,j)$ for attacker $A$ if the edge is or is not interdicted, respectively. We assume that the realizations of these belief probabilities are independent of each other. We further assume that the probabilities $\delta_D(i,j,A)$ and $\tilde{\delta}_D(i,j,A)$ are independent between edges.


To derive the \textit{OptDef} formulation, note that a level $k$ defender's portfolio interdicts a set of edges $S\subseteq E$. Further, let $P(i)$ be the set of paths leading from the source node $1$ to node $i>1$. Then, the conditional probability of a successful attack is
\begin{align}
    \label{cyber:eqn:def_path_expect}
    &   
    \sum_{A \in \mathcal{A}}  \theta_D(A) \sum_{p \in P(n)} \psi^k_D(p,A)  \prod_{(i,j)\in p\setminus S}\delta_D(i, j, A) \prod_{(i,j)\in p\cap S} \tilde{\delta}_D(i, j, A).
\end{align}
The maximum probability path from any node $i\in N$ to the sink node $n$  does not depend on the path used to reach $i$.  
Therefore, can rewrite (\ref{cyber:eqn:def_path_expect}) as
\begin{align}
    \label{cyber:eqn:X0} 
    \sum_{A \in \mathcal{A}} \theta_D(A)  \sum_{p \in P(n)} \prod_{(i,j)\in p}\psi^k_D(i, j, A) \left( \delta_D(i, j, A) (1 - x_{ij} ) + \tilde{\delta}_D(i, j, A)x_{ij} \right)
\end{align}

The defender believes that the probability that attacker $A$ attempts to traverse edge $(j,\ell) \in E$ is $y_{ij}(A)$ if the edge is not interdicted or $\tilde{y}_{j \ell}(A)$ if it is interdicted, which are decision variables in the defender's model. We define $y_{j \ell}(A)$ and $\tilde{y}_{j \ell}$(A) recursively as
\begin{align*}
    y_{j \ell}(A)
    & = (1 - x_{j \ell})\psi^k_D(j, \ell, A) \sum_{(i,j) \in E^-_j} \left( \delta_D(i,j,A)y_{ij}(A) + \tilde{\delta}_D(i,j,A) \tilde{y}_{ij}(A)\right)\\
    \tilde{y}_{j \ell}(A)  
    & = x_{j \ell}\psi^k_D(j, \ell, A) \sum_{(i,j) \in E^-_j} \left( \delta_D(i,j,A)y_{ij}(A) + \tilde{\delta}_D(i,j,A) \tilde{y}_{ij}(A)\right) .
\end{align*}



Gathering these equations, we introduce the \textit{OptDef} formulation:
\begin{align}
z_D^{k} = 1 - \min & \sum_{A\in \mathcal{A}}\theta_D(A) \sum_{(i,n) \in E^-_n}\left( \delta_D(i,n,A)y_{in}(A) + \tilde{\delta}_D(i,n,A) \tilde{y}_{in}(A)\right) & \label{cyber:eqn:onj}
\end{align}
\vspace{-30pt}
\begin{align}
\text{s.t.  }& \label{cyber:eqn:D1} y_{1i}(A)+\tilde{y}_{1i}(A) \geq \psi^k_D(1,i,A), & \forall A\in \mathcal{A}, (1, i) \in E^+_1\\
& \nonumber y_{j \ell}(A)+\tilde{y}_{j \ell}(A)  \geq \psi^k_D(j,\ell,A) \sum_{(i,j) \in E^-_j}\left( \delta_D(i,j,A)y_{ij}(A) + \tilde{\delta}_D(i,j,A) \tilde{y}_{ij}(A)\right),\hspace{-80mm} & \\[-4\jot]
& & \hspace{-80mm} \forall A\in \mathcal{A}, j \in N\setminus\{1,n\}, (j,\ell) \in E^+_j \label{cyber:eqn:D2}
\end{align}
\vspace{-40pt}
\begin{align}
& \label{cyber:eqn:D4} \tilde{y}_{ij}(A) \leq \sum_{m\in M_{ij}}w_m, &\forall A\in \mathcal{A}, (i,j) \in E\\
& \label{cyber:eqn:D5} \sum_{m\in M}\kappa_D(m) w_m \leq \beta_D & \\
& \label{cyber:eqn:D6} y_{ij}(A),\tilde{y}_{ij}(A) \ge 0, &\forall A\in \mathcal{A},\ (i,j) \in E \\
& \label{cyber:eqn:D7} w_{m} \in \{0,1\}, &\forall m\in M 
\end{align}

A level $k$ defender maximizes the conditional probability that an attack is prevented (\ref{cyber:eqn:onj}), which is equivalent to minimizing the conditional probability that an attack succeeds, and which occurs when an attacker traverses the network. 
Constraint set (\ref{cyber:eqn:D1}) enforces that each attacker attempts an attack. Constraint set (\ref{cyber:eqn:D2}) serves as a flow balance equation, determining the probability that an attack from each attacker reaches each edge. Specifically, this is the probability of an attack reaching the edge's starting node multiplied by the conditional probability of the attack then progressing on that edge. Constraint set (\ref{cyber:eqn:D4}) allows edges to be interdicted only if a control is chosen that interdicts that edge.  Constraint (\ref{cyber:eqn:D5}) allows the defender to choose only as many controls as their budget allows. The last two sets of constraints, (\ref{cyber:eqn:D6}) and (\ref{cyber:eqn:D7}), ensure that the appropriate variables are non-negative and binary.


A solution to the level $k$ defender's problem informs the level $k+1$ attackers' beliefs regarding the defender's edge reliabilities. In particular, $\Omega_A^{k}$ is defined as the set $\{m \in M \; : \; w_m = 1\}$ and for each $(i,j) \in E$  
\begin{align}\label{cyber:eqn:Dupdate}
    \delta_A(i,j|\Omega^k_A)=\begin{cases}
      \delta_D(i,j,A), & \text{if}\ x_{ij} = 0 \\
      \tilde{\delta}_D(i,j,A), & \text{if}\ x_{ij} = 1.
    \end{cases}
\end{align}

%
\section{Approximation of the Defender's Problem}\label{sec:algorithms}

It may be computationally difficult to find an optimal solution to \textit{OptDef} for large-scale problem instances. However, we show that there exists an approximation algorithm that provides a solution with a $1-1/e$ performance guarantee. 
To do so, we show that \textit{OptDef} is equivalent to a non-negative, submodular maximization problem with a knapsack constraint.  A heuristic can find solutions that are at least $(1-1/e)$ of the optimal solution in polynomial time for this problem \citep{sviridenko_note_2004}.  A similar result exists for the related problem of maximizing a submodular function subject to a cardinality constraint \citep{nemhauser_analysis_1978}.
Later in Section \ref{cyber:sec:computational}, we show that this approximation algorithm often identifies solutions whose values are extremely very close to the optimal solution value in practice. This approximation algorithm may be used to decrease the time required to find an optimal solution to \textit{OptDef} by providing a warm start for the MIP solver.

We adapt the heuristic from \citet{khuller1999budgeted} to \textit{OptDef}. We begin by selecting the maximum value set of controls, $S_1 \subset M$, with $| S_1 | \le 2$. Next, we enumerate all sets of controls with a cardinality of three, $\mathcal{S_2}$, that satisfy the budgetary constraint, i.e., we have
$\sum_{m\in S_2}\kappa_D(m) \leq \beta_D$ for each set of controls $S_2 \in \mathcal{S_2}$ with $| S_2 | = 3$. Then, we greedily complete each of these sets until the budget or available controls are exhausted.  That is, for a submodular function $f$ on a set $M$, we add an element $m^* \in M \setminus S_2$ if it satisfies the budget constraint (i.e., $\kappa_D(m^*) + \sum_{m\in S_2}\kappa_D(m) \leq \beta_D$) that satisfies
\begin{align}\label{greedyrule}
    m^* \in \argmax_{m\in M\setminus S_2}\frac{f\left(S_2\cup\{m\}\right)-f\left(S_2 \right)} {\kappa_D(m)}.
\end{align}
Then, $S_2 \leftarrow S_2 \cup m^*$ and the process is repeated until no new elements can be selected. Let the maximum value set completed this way be $S_3$.  If $f(S_1) \geq f(S_3),$ then the algorithm returns $S_1.$ Otherwise, it returns $S_3.$

Theorem \ref{cyber:lem:supermod} demonstrates that the conditional probability of a successful attack, the complement of the \textit{OptDef} objective function value, is a supermodular function given a set of interdicted edges. 

\begin{thm}
\label{cyber:lem:supermod} 
Let $S$ be a set of interdicted edges $(i,j) \in E$ and let 
\begin{align*}
    0\leq \tilde{\delta}_D(i, j, A) \leq \delta_D(i,j,A), \forall (i, j)\in E, A\in \mathcal{A}.
\end{align*}
Then for any path $p\in P$ and attacker $A\in \mathcal{A},$ 
\begin{align*}
    g(S, p, A) &= \prod_{(i,j)\in p\setminus S}\delta_D(i, j, A) \prod_{(i,j)\in p\cap S} \tilde{\delta}_D(i, j, A)
\end{align*} 
is a non-increasing, supermodular function in $S$.
\end{thm}
\proof
We begin by showing that $g(S,p,A)$ is non-increasing in $S$. First, note that $0\le g(S,p,A) \le 1$, since $g(S,p,A)$ is a product of numbers whose values are between $0$ and $1$. For sets of edges $S_1\subseteq S_2\subseteq E$, we have
\begin{equation}
\begin{aligned}[b]
     g(S_2,p,A) & = \prod_{(i,j)\in p\setminus S_2}\delta_D(i, j, A) \prod_{(i,j)\in p\cap S_2} \tilde{\delta}_D(i, j, A)\\
     & = \prod_{(i,j)\in p\cap(S_2\setminus S_1)}\frac{\delta_D(i, j, A)\tilde{\delta}_D(i, j, A)}{\delta_D(i, j, A)\tilde{\delta}_D(i, j, A)}\prod_{(i,j)\in p\setminus S_2}\delta_D(i, j, A) \prod_{(i,j)\in p\cap S_2} \tilde{\delta}_D(i, j, A)\\
    & = \prod_{(i,j)\in p\cap(S_2\setminus S_1)}\frac{\tilde{\delta}_D(i, j, A)}{\delta_D(i, j, A)}\prod_{(i,j)\in p\setminus S_1}\delta_D(i, j, A) \prod_{(i,j)\in p\cap S_1} \tilde{\delta}_D(i, j, A)\\
    & = \left(\prod_{(i,j)\in p\cap(S_2\setminus S_1)}\frac{\tilde{\delta}_D(i, j, A)}{\delta_D(i, j, A)}\right)g(S_1,p,A). \label{cyber:eqn:non-inc} 
\end{aligned}
\end{equation}
Since $0\leq \tilde{\delta}_D(i, j, A)\leq \delta_D(i,j,A), \forall (i, j)\in E$, (\ref{cyber:eqn:non-inc}) shows that $g(S,p,A)$ is non-increasing in $S$.

We now prove that $g(S,p,A)$ is supermodular by showing for all $(i_1,j_1),(i_2,j_2) \in E$ that 
\begin{align}
    \nonumber&g(S, p, A) + g(S \cup \{(i_1,j_1), (i_2,j_2)\}, p, A)\\ - 
    \label{cyber:eqn:supermod_defn}&g(S \cup \{(i_1,j_1)\}, p, A) - g(S \cup \{(i_2,j_2)\}, p, A)\geq 0.
\end{align}

We use (\ref{cyber:eqn:non-inc}) to algebraically simplify the left-hand side of (\ref{cyber:eqn:supermod_defn}):
\begin{align}
   \nonumber &g(S, p, A) + g(S \cup \{(i_1,j_1), (i_2,j_2)\}, p, A) - 
    g(S \cup \{(i_1,j_1)\}, p, A) - g(S \cup \{(i_2,j_2)\}, p, A)\\
   \nonumber =& g(S, p, A) \left(1 + \frac{\tilde{\delta}_D(i_1, j_1, A)\tilde{\delta}_D(i_2, j_2, A)}{\delta_D(i_1, j_1, A)\delta_D(i_2, j_2, A)} -\frac{\tilde{\delta}_D(i_1, j_1, A)}{\delta_D(i_1, j_1, A)} - \frac{\tilde{\delta}_D(i_2, j_2, A)}{\delta_D(i_2, j_2, A)} \right)\\
  \label{cyber:eqn:supermod_simple}  = & g(S, p, A) \left(1 - \frac{\tilde{\delta}_D(i_2, j_2, A)}{\delta_D(i_2, j_2, A)} \right) \left(1-\frac{\tilde{\delta}_D(i_1, j_1, A)}{\delta_D(i_1, j_1, A)}\right).
\end{align}

Since, $0\leq \tilde{\delta}_D(i, j, A)\leq \delta_D(i,j,A), \forall (i, j)\in E$, we have the following three inequalities:
\begin{align*}
g(S, p, A) & \geq 0,\\
1 - \frac{\tilde{\delta}_D(i_2, j_2, A)}{\delta_D(i_2, j_2, A)} & \geq 0,\\
1 - \frac{\tilde{\delta}_D(i_1, j_1, A)}{\delta_D(i_1, j_1, A)}  & \geq 0.
\end{align*}
Combining these three inequalities with (\ref{cyber:eqn:supermod_simple}) yields (\ref{cyber:eqn:supermod_defn}).  Therefore, $g(S, p, A)$ is a non-increasing supermodular function. $\blacksquare$


The consequence of Theorem 1 is that the \textit{OptDef} objective function is a non-negative, submodular function. Next, we prove that \textit{OptDef} is a non-negative, submodular maximization problem subject to a knapsack constraint.

\begin{thm}\label{cyber:thm:supermod_optdef} \textit{OptDef} is a non-negative, submodular maximization problem subject to a knapsack constraint.
\end{thm} 
\proof  
We begin by reformulating (\ref{cyber:eqn:X0}), which is equivalent to the complement of the \textit{OptDef} objective function the objective (\ref{cyber:eqn:onj}) and Constraints (\ref{cyber:eqn:D1}) -- (\ref{cyber:eqn:D2}). Then the probability that an attack is prevented, using the definition of $g(S,p,A)$ from Theorem \ref{cyber:lem:supermod} is
\begin{align}\label{thmeqn}
     1 - \sum_{A \in \mathcal{A}} \theta_D(A)  \sum_{p \in P(n)} g(S,p,A)\prod_{(i,j)\in p}\psi_D(i, j, A).
\end{align}
where $S$ is the set of edges interdicted by the defender's portfolio. Note that $\prod_{(i,j)\in p}\psi_D(i, j, A)$ is a constant in this formulation. The non-negative weighted sum of supermodular functions is also supermodular, yielding a function in \eqref{thmeqn} that is submodular after multiplying by $-1$ and subtracting the term from 1. The resulting value is between $0$ and $1$, since it represents a probability.

The remaining constraints of \textit{OptDef}, Constraints (\ref{cyber:eqn:D4}) -- (\ref{cyber:eqn:D5}), simply define a knapsack constraint on the set of controls. Since $S$ is non-decreasing in the controls, this is equivalent to a knapsack constraint on $S$. Hence, we have the desired result. $\blacksquare$

There are two implications of Theorem 2. The first is that the algorithm introduced earlier in this section from \cite{khuller1999budgeted} identifies solutions with a guaranteed $1-1/e$ approximation ratio. 
The second implication is that a more computationally efficient greedy algorithm (also introduced by \cite{khuller1999budgeted}) can be used to identify solutions with an approximation ratio of $1-1/\sqrt{e}$. 
The greedy algorithm starts with an empty set of controls and greedily completes this set according to \eqref{greedyrule} until the budget or available controls are exhausted. Then it compares this greedily chosen set with the maximum value single-element set and chooses the one with the higher objective function value.


\section{Iterative ARA Algorithm}
\label{sec:araFramework}
We formally introduce an iterative algorithm for solving \textit{OptAtt} and \textit{OptDef} for the attacker and defender problems across all values of $k$.

\begin{enumerate}
\item Initialization: Provide data.

Some sets and beliefs given in Table \ref{tab:systemnotation} must be provided, including the level ($K$) of defender, the set of attackers ($\mathcal{A}$), sets of nodes ($N$) and edges ($E$), set of controls ($M$), reliabilities of edges ($\delta_D(i,j,A)$, $\Tilde{\delta}_D(i,j,A) \;\forall (i,j) \in E, \,A \in \mathcal{A}$), defender budget ($\beta_D$), costs of controls ($\kappa_D(m) \; \forall \,m\in M$), and the distribution of attackers ($\theta_D(A) \;\forall\,A\in\mathcal{A}$).

\item Identify level 0 attacker and defender solutions:

We assume a level 0 attacker $A$ takes a greedily formed path $p$ as described in Section \ref{cyber:sec:optatt}.
This yields $\psi_D^{k'}(i,j,A) = 1$ if $(i,j) \in p$ and $\psi_D^{k'}(i,j,A) = 0$ if $(i,j) \notin p$ for all defenders $k' > 0$ who believe that attacker $A$ is at level 0. 

We assume a level 0 defender does not choose any controls, so that the level 0 defender's solution is $\Omega_A^0 = \emptyset$ with $\delta_A^1(i,j|\Omega_A^0) = \delta_D(i,j,A)$ for each attacker $A \in \mathcal{A}$. 

\item For $k = 1,2,...,K-1$:

\begin{itemize}
\item Solve \textit{OptAtt} using Dijkstra's algorithm:

For each level $k$ attacker $A$ we solve \textit{OptAtt}, given by \eqref{cyber:eqn:A0} -- \eqref{cyber:eqn:A7}, to obtain the attack path $p$ using Dijkstra's algorithm as described in Section \ref{cyber:sec:optatt}. Use the values of $u_{ij}$ to set $\psi_D^{k'}(i,j,A), \ (i,j) \in E$ for each $k<k'\leq K$ using \eqref{eq:defenderupdate}.

\item Solve \textit{OptDef} for a level $k$ defender. We suggest two ways to do this:
\begin{itemize}
    \item Solve \textit{OptDef} to optimality using the mixed integer programming formulation given by (\ref{cyber:eqn:onj}) -- (\ref{cyber:eqn:D7}).
    \item Identify approximate solutions using the approximation algorithm or the greedy algorithm presented in Section \ref{sec:algorithms} to identify solutions within $(1-1/e)$ or $(1-1/\sqrt{e})$ of the optimal solution value, respectively.
\end{itemize}
Given a solution to \textit{OptDef}, let $\Omega_A^k = \{m \in M \;|\;w_m= 1\}$. Then set the values of $\delta_A^{k+1}(i,j|\Omega^k_A)$ for all level $k+1$ attackers $A \in \mathcal{A}$ using \eqref{cyber:eqn:Dupdate}. 
\end{itemize}
	
\item Solve \textit{OptDef} a final time to obtain the solution for a level $K$ defender.
\end{enumerate}

\added{To illustrate the algorithm, we revisit the example introduced in Section \ref{subsec:approach}. Again, we assume the defenders budget allows them to select only a single control.
Next, we iterate through the algorithm to find the OptDef solutions with level 1 -- 4 defenders, respectively. For each level $k$ defender from $k=1$ to $k=4$, we assume the defender is defending against a uniform distribution of the attackers of levels less than $k$. We report the attack success probability for both models to directly compare the \textit{OptAtt} and \textit{OptDef} solutions in Table \ref{tab:example}. }

\added{We first consider level 0 attackers, who choose the right path using a greedy tactic, with an attack success probability of 0.102. A level 0 defender selects no controls. A level 1 attacker optimizes against the level 0 defender and again chooses the right path for an attack success probability of 0.102, since this is the maximum reliability path with all uninterdicted edges.}

\added{We continue to consider each level of attacker and defender up to the desired defender level. Next, we solve \textit{OptDef} for a level 1 defender, who optimizes against level 0 attackers. The level 1 defender chooses the control that most decreases the reliability of the level 0 attackers' path, which is control $m_3$. This attains an \textit{OptDef} attack success probability of 0.074.
The level 2 attacker then selects the left path, yielding an attack success probability of 0.09.
A level 2 defender optimizes against both level 0 and level 1 attackers, who both take the right path, and again selects $m_3$, giving an attack success probability of 0.074.
Given this, a level 3 attacker takes the right path with an attack success probability of 0.102.}

\added{Next, we consider a level 3 defender, who defends against level 0 -- 2 attackers, two of whom take the right path and one of whom takes the left path. A level 3 defender selects $m_1$ for an \textit{OptDef} attack success probability of 0.078. The level 4 defender defends against all levels 0 -- 3 attackers, two of whom traverse the left path and two of whom traverse the right path. 
In this case, $m_2$ becomes the optimal control choice. This achieves a final \textit{OptDef} attack success probability of 0.069.}
In this example we see the defender change their control decisions based on what paths they believe the attackers would take. 

\begin{table}[htp]\centering \caption{Example \textit{OptAtt} and \textit{OptDef} solutions for attackers and defenders of varying levels}\label{tab:example}
\begin{tabular}{|ccc|ccc|}
\hline \multicolumn{3}{|c}{Attacker} & \multicolumn{3}{|c|}{Defender} \\ \hline
 & Path & \textit{OptAtt} Attack & & Mitigation & \textit{OptDef} Attack  \\ 
k & selected & Success Probability & k & selected & Success Probability \\ \hline
0 & Right & 0.102 & 0 & N/A   & N/A\\
1 & Right & 0.102 & 1 & $m_3$ & 0.074\\
2 & Left  & 0.09  & 2 & $m_3$ & 0.074 \\
3 & Right & 0.102 & 3 & $m_1$ & 0.078\\
--&--     &--     & 4 & $m_2$ & 0.069\\ 
\hline 
\end{tabular}
\end{table}

Note that the ARA algorithm results in a series of solutions to \textit{OptDef} that correspond to level $0, 1,...,K$ defenders. Therefore, our approach provides decision-makers with a suite of solutions instead of a single, ``good'' solution, which could be advantageous for defensive planning.

\section{Computational Results}
\label{cyber:sec:computational}




In this section, we provide computational results for the models and algorithms based on the information security investment case study introduced by \cite{dubois2020optimizations}. 
In the case study, we generate random data according to a specific parameterized structure that we are able to adjust to leverage control over the size and complexity of the problem. We describe how this structure is created as well as the parameters used. Later, we vary the size of the problem instances to assess the performance of the approximation algorithm.

\begin{figure}[ht!]
    \centering
    \includegraphics[width=8cm]{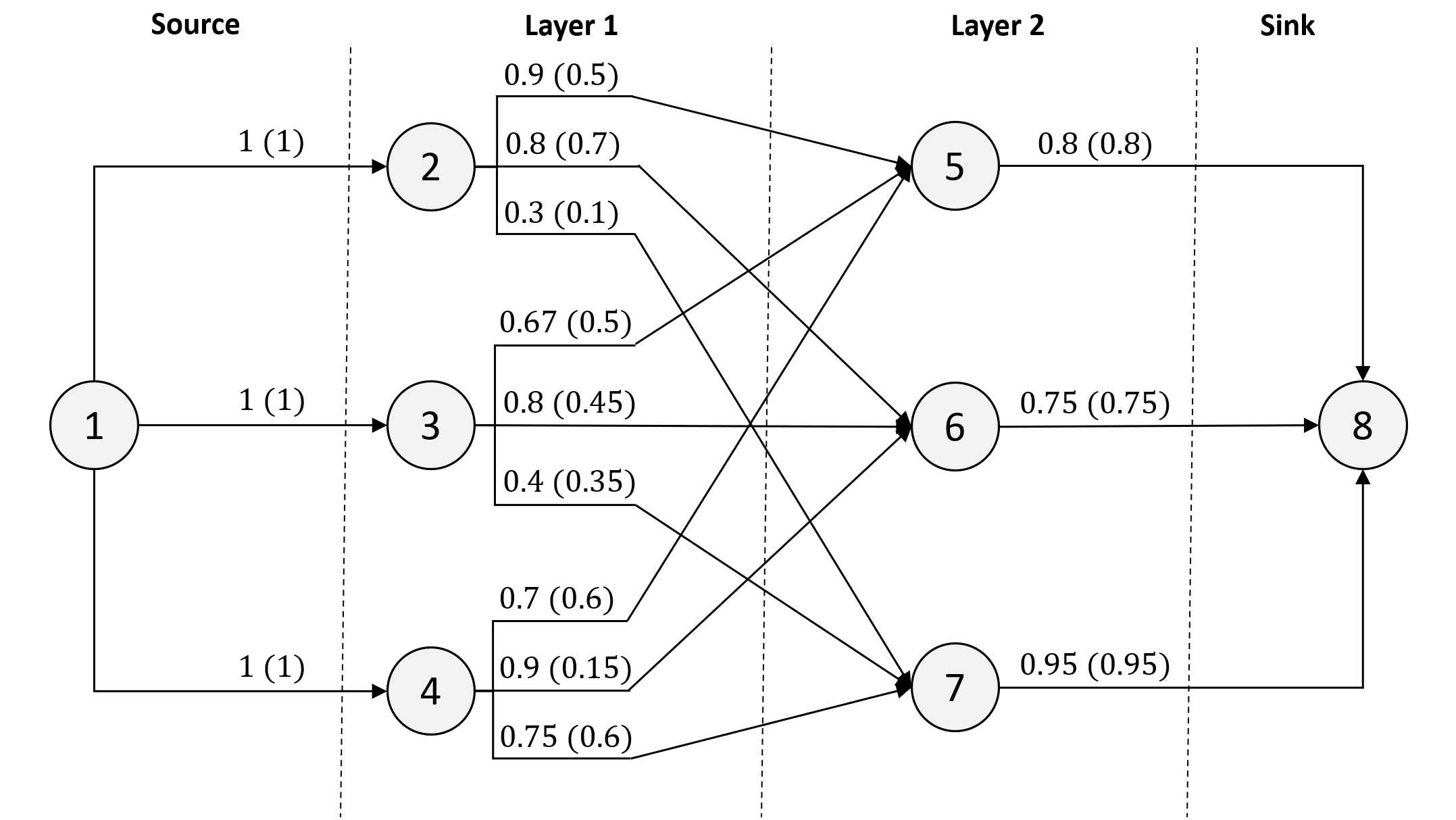}
    \caption{Simple example of synthetic network structure with two layers.  Example edge reliabilities are shown along each edge for when the edge is not interdicted (interdicted) , $\delta_{ij}(i,j,A) \left(\tilde{\delta}_{ij}(i,j,A)\right)$}
    \label{cyber:fig:Synthetic}
\end{figure}

In the case study, each problem instance is given by a network that is organized into a source node, a sink node, and $\ell$ ``layers'' of nodes.  Each layer of nodes is a set of nodes such that nodes in a layer connect only to nodes in the next layer. For simplicity, we assume that each layer contains the same number of nodes, which results in a structure with many interconnections and paths. 
Figure \ref{cyber:fig:Synthetic} illustrates an example of the network under consideration with two layers of three nodes each. A path through the network includes one node from each layer.  The source node is connected by non-interdictable edges to each node in the first set, with reliabilities of $1$.  Three edges leave each node in sets $1,..., \ell-1$ to nodes in sets $2,..., \ell$, respectively. 

We first assign edges such that each node has both an entering and leaving edge, before randomly assigning the remaining edges. We then assign reliabilities to each of these edges randomly, and each reliability is treated as a constant once it is generated. The un-interdicted reliabilities are uniformly assigned a value between 0 and 1, i.e., $\delta(i,j)\sim U(0,1]$.  The interdicted reliabilities are a random $(0,1)$ proportion of the un-interdicted value, i.e., $\tilde{\delta}(i,j) \sim U(0,1]\cdot \delta(i,j)$.  Finally, non-interdictible edges connect each node in layer $\ell$ to the sink node, with reliabilities of $U(0.5,1]$. We choose a number of possible controls and a budget for the defender, then randomly determine costs for each of these controls using a uniform distribution over some specified cost range.  This cost is exactly $1$ for a cardinality budget constraint and the costs follow a $U[0.5,1.5]$ distribution for a knapsack constraint to maintain an average cost of $1.0$ per control in both situations. Each control is then randomly assigned a subset of these edges to interdict. To do so, we first specify a parameter, $\alpha$, that gives the average proportion of edges covered by any control. If we have a cardinality budget constraint, for each control we loop through the edges and assign each to the control with a probability of $\alpha$. If we have a knapsack constraint, we alter the probability of assigning each edge to a control by considering the cost. Let $\Delta_m$ be the value attained by subtracting the average cost from the cost of the control $m$, then dividing by the range of possible costs. In our case with costs generated from the range $[0.5,1.5]$, this is equivalent to $\Delta_m = \kappa_D(m) - 1$. Then we assign edges to control $m$ with a probability of $\alpha(1 + \alpha_2\Delta_m)$, where $\alpha_2$ represents the amount of effect cost has on mitigation quality.

For our case study we generated a problem instance for a graph with $\ell =5$ layers and 5 nodes per layer. To better visualize the tradeoffs between controls, a cardinality constraint was used for the defender's budget, where the defender can choose 4 of 10  controls. Each control was randomly generated to cover approximately 15\% of the edges ($\alpha=0.15$). \added{While a realistic value of the defender's level is 4, we set a maximum defender level of $K = 10$ to study a range of defensive solutions}. Each defender level $k \leq K$ assumes a uniform distribution of attacker levels over $0 \le k'<k$. 

We solve the case study using the algorithm in Section \ref{sec:araFramework}, where for each level of defender we solve \textit{OptDef} to optimality using Gurobi 9.1.1. This computation was performed on a computer with an Intel(R) Core(TM) i7-8650U CPU @1.90GHz  2.11GHz processor and 16GB of RAM, and took 1.36 seconds to run the full algorithm.

The ARA algorithm solves the defender's problem eleven total times across all levels of $k$, thereby generating various potential defender solutions before arriving at the final solution for the level $K = 10$ defender. 
Table \ref{tab:defSolns} reports the portfolio of controls chosen by each defender level $k$ from 0 to 10. Note that by assumption the level 0 defender chooses no controls. 
%
%
Table \ref{tab:defSolns} indicates that some controls are not chosen in any of the defender's solutions, such as $m_1$, $m_6$, $m_7$, and $m_8$. Other controls work together to provide better protection against different or more of a variety of paths the attackers may take (e.g., controls $m_2$, $m_3$, and $m_4$). Controls $m_5$ are $m_9$ are chosen by defenders whose level is at least 3 and 5, respectively, indicating that some controls only become attractive by more strategic defenders (who believe they face more strategic attackers).

\begin{table}[h]
\caption{Table of portfolios selected by the defender for differing levels of $k$}\centering
    \begin{tabular}{c|llllllllll}
        & \multicolumn{9}{c}{Controls} \\ \hline
        Level of Defender & $m_1$ & $m_2$ & $m_3$ & $m_4$ & $m_5$ & $m_6$ & $m_7$ & $m_8$ & $m_9$ & $m_{10}$ \\ \hline
        0                   & - & - & - & - & - & - & - & - & - & - \\
        1                   & - & \checkmark & \checkmark & \checkmark & - & - & - & - & - & \checkmark \\
        2                   & - & \checkmark & \checkmark & \checkmark & - & - & - & - & - & \checkmark \\
        3                   & - & \checkmark & - & \checkmark & \checkmark & - & - & - & - & \checkmark \\
        4                   & - & \checkmark & - & \checkmark & \checkmark & - & - & - & - & \checkmark \\
        5                   & - & \checkmark & - & \checkmark & \checkmark & - & - & - & \checkmark & - \\
        6                   & - & \checkmark & - & \checkmark & \checkmark & - & - & - & \checkmark & - \\
        7                   & - & \checkmark & \checkmark & - & \checkmark & - & - & - & \checkmark & - \\
        8                   & - & \checkmark & \checkmark & - & \checkmark & - & - & - & \checkmark & - \\
        9                   & - & - & \checkmark & \checkmark & \checkmark & - & - & - & \checkmark & - \\
        10                  & - & - & \checkmark & \checkmark & \checkmark & - & - & - & \checkmark & -
    \end{tabular}
    \label{tab:defSolns}
\end{table}



\begin{figure}[ht!]
    \centering
    \begin{subfigure}{8cm}
        \includegraphics[width=8cm]{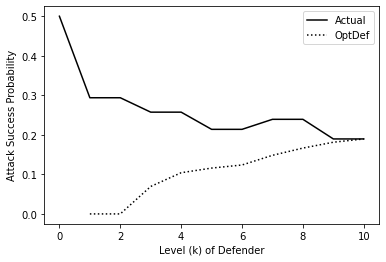}
        \caption{Example Case}
    \end{subfigure}
    \begin{subfigure}{8cm}
        \includegraphics[width=8cm]{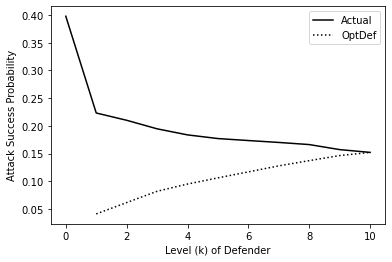}
        \caption{Averaged over 100 instances}
    \end{subfigure}
    \caption{Two lines showing the average attack success of attackers levels 0-9 against levels $k < 10$ defenders, and the OptDef value for levels $k< 10$ defenders (attack success of a uniform distribution of attackers levels $0$ to $k$).}
    \label{fig:optdefVsActual}
\end{figure}

Next, we examine the \emph{attack success probability},  the average probability of a successful attack and the complement to the \textit{OptDef} objective function value. 
We report the attack success probability for defenders with levels $0-9$ given a uniform distribution of levels $0 - 9$ attackers. Figure \ref{fig:optdefVsActual} reports the attack success probability from both the defender's and attackers' points of view. First, it reports the attack success probability corresponding to \textit{OptDef} (the dotted line), which captures the attack success probability based on the defender's belief regarding the attackers. Second, it reports the actual, retrospective attack success probability for the same set of attackers based on their actual levels of strategic sophistication (the solid line). These values differ, since defenders with levels lower than 10 are incapable of perceiving level 9 (or higher) attackers' true levels of strategic sophistication.

Figure \ref{fig:optdefVsActual}(a) illustrates the results across a single problem instance, and Figure \ref{fig:optdefVsActual}(b) illustrates the results averaged over 100 problem instances of the same size as the example problem. 
The attack success probability associated with \textit{OptDef} is increasing in $k$ in both figures. This occurs since more strategic defenders select controls for more strategic attackers who are better able to evade defenses. 
%
We observe that the actual attack success probability is higher than that associated with \textit{OptDef}, since defenders with levels $k<10$ believe that some attackers are less strategic than they are in actuality, resulting in 
an inaccurate believed attack success probability. Both lines converge to the same point when the defender's beliefs match reality, which in our example occurs when the defender is at level $k=10$.
Figure \ref{fig:optdefVsActual}(a) indicates that the actual attack success probability does not always decrease with $k$ for any particular problem instance. This effect is due to attackers with higher levels of $k$. Since the defender does not accurately assume the distribution of attackers, it is possible that attackers can find more reliable paths against a sophisticated defender.

Figure \ref{fig:optdefVsActual}(b) illustrates the attack success probability associated with \textit{OptDef} and for the actual values averaged across 100 randomly generated problem instances of the same size.
The overall trends in Figure \ref{fig:optdefVsActual}(b) are similar to those in Figure \ref{fig:optdefVsActual}(a), with Figure \ref{fig:optdefVsActual}(b) showing that the average actual attack success probability monotonically decreases with $k$. This suggests that the solutions provided by \textit{OptDef}, on average, provide better defenses for more strategic defenders.



Next, we consider the effects of the defender's portfolio of controls and how it performs against attackers with different levels of strategic sophistication. As before, we consider defenders of levels 0 through 9. We study how two subsets groups of attackers perform against these defenders: low-level (levels $0-4$) and high-level (levels $5-9$). We also report all attackers (levels $0-9$) for comparison. Note that the case with level $0-9$ attackers matches that considered in Figure \ref{fig:optdefVsActual}.
In each case, the attackers' actual levels follow a discrete uniform distribution as before. Each defender assumes that the levels of the attackers are uniformly distributed across all levels $k$ lower than that of the defender. 
The level 5 defender believes they are defending against low-level $0-4$ attackers, and a level 10 defender believes they are defending against level $0-9$ attackers. The defender's beliefs matches the actual attacker levels in both of these cases, however, there is a mismatch between the defender's beliefs and the actual levels of the attackers in all other cases, where the defender either overestimates or underestimates the attackers' levels to varying degrees.

\begin{figure}[ht!]
    \centering
    \begin{subfigure}{8cm}
        \includegraphics[width=8cm]{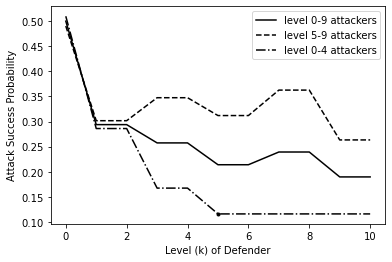}
        \caption{Example Problem}
    \end{subfigure}
    \begin{subfigure}{8cm}
        \includegraphics[width=8cm]{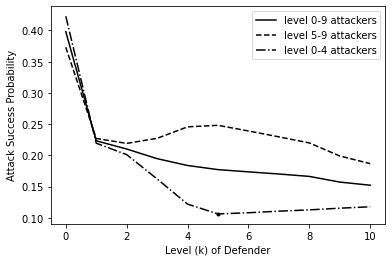}
        \caption{Averaged over 100 instances}
    \end{subfigure}
    \caption{Attacker success probabilities against defenders of levels 0 through 9, for three uniform distributions of attackers: low-level (0-4), high-level (5-9), and both (0-9).}
    \label{fig:exAttackSuccessByLevel}
\end{figure}

Figure \ref{fig:exAttackSuccessByLevel} illustrates the actual attack success probability for the defender across the three groups of attackers, with Figure \ref{fig:exAttackSuccessByLevel}(a) illustrating the actual attack success probability for the example problem and Figure \ref{fig:exAttackSuccessByLevel}(b) illustrating the actual attack success probability averaged across 100 problem instances. In both figures, higher level attackers tend to be more successful than lower level attackers across all defender levels, which is expected. We observe that that higher level defenders tend to achieve a lower attack success probability when facing level $0-4$ attackers, which suggests that defenses may be more effective against less strategic attackers, 
even when the defender overestimates the attackers' levels (which occurs for level $6-9$ defenders).  

In contrast, level $5-9$ attackers often achieve high attack success probabilities. In Figure \ref{fig:exAttackSuccessByLevel}(a), the attack success probability of level $5-9$ attackers is highest against level $7$ and $8$ defenders. This is surprising and occurs since the attackers whose levels are higher than the defender can find paths that the defender does not defend well. This indicates that the defenses chosen by more sophisticated defenders may sometimes perform worse against sophisticated attackers than the defenses chosen by less sophisticated defenders.

Figure \ref{fig:exAttackSuccessByLevel}(b) presents results averaged across 100 randomly generated problem instances. Overall, we see a decreasing trend in the attack success probability across all level $0-9$ attackers and low-level $0-4$ attackers as the defender's level increases.
The decreasing attack success probability is more pronounced for level $0-4$ attackers. This occurs, since the defender increasingly assumes a more accurate distribution of the attackers as the defender level increases to $5$. After level 5, the defender continues to defend well across level $0-4$ attackers and mostly improves against level $5-9$ attackers. 
The actual attack success probability against the high-level $5-9$ attackers does not monotonically decrease with the defender's level (i.e., it increases between defender levels $2$ to $5$). This is counter-intuitive, and results from the defender's beliefs. A level $k=5$ defender, for example, assumes they face level $0-4$ attackers and therefore optimizes defenses for less strategic attackers. As the defender's level increases to $k=9$, the defender becomes more effective against the most strategic attackers. Overall, this example suggests that is better for the defender to overestimate the level of the attackers rather than to underestimate the level of the attackers.

\begin{figure}
    \centering
    \includegraphics[width=10cm]{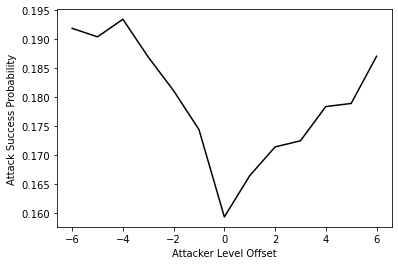}
    \caption{Attack success of a uniform distribution of attackers against level 10 defenders who believe all attackers' levels are offset by the amount given on the x-axis.}
    \label{fig:overunderestimation}
\end{figure}

To further analyze the effects of the defender misjudging the sophistication of the attackers, we consider a defender of level $k=10$ who uniformly overestimates or underestimates the attackers' levels by a certain \emph{offset}. As before, we consider a uniform distribution of attackers whose true levels are 0 through 9. Positive offsets correspond to the defender overestimating the attackers' levels, and negative offsets correspond to underestimating the attackers' levels. When overestimating, the defender cannot perceive attackers of a higher level than them, and thus they assume any such attackers are one level below them in \textit{OptDef}. Likewise, when underestimating, the defender always assumes the attackers have a level of at least 0. A defender with an offset of -3, for example, only perceives level 0 through 6 attackers after capping the perceived attacker levels below at 0.

Figure \ref{fig:overunderestimation} depicts the actual attack success probability averaged across 100 problem instances as a function of the offset. We observe that both over and underestimating the distribution of the attackers leads to higher attack success probabilities in comparison to the case when the defender perceives the true attacker levels (shown by the offset of 0), with underestimating the attackers being slightly worse than overestimating. 
This observation is consistent with Figure \ref{fig:exAttackSuccessByLevel}. This suggests that it is ideal for the defender to ``get it right'' and correctly model attackers and that there is a benefit to erring on the side of planning for sophisticated attackers.

We next assess the quality of the solutions identified by the greedy algorithm (see Section \ref{sec:algorithms}), which is guaranteed to identify solutions whose objective function values are at least $1-1 / \sqrt{e}$ of the optimal \textit{OptDef} solution values.
We consider randomly generated instances for a level 10 defender, each using a knapsack constraint for the budget. 
We compare the greedy and exact algorithms on the same problem instances with a level $10$ defender as follows by first using the greedy algorithm to identify near-optimal \textit{OptDef} solutions for defenders of levels $1-9$. Then, we use either an exact algorithm or the greedy algorithm for the level 10 defender's problem instance.

We randomly generated $40$ problem instances with varying sizes to evaluate the solution quality of the greedy algorithm. Table \ref{tab:algm} reports the parameters used to generate each problem instance as well as the CPU times associated with the exact (\textit{OptDef}) and greedy algorithms. The CPU time (in seconds) to solve \textit{OptDef} to optimality includes the setup time of writing the variables and constraints as well as the solver time. Likewise, the time to implement the greedy algorithm includes the time to read the inputs and execute the algorithm. Table \ref{tab:algm} reports the ratio of the greedy solution value ($z_h$) and the optimal solution value to \textit{OptDef} ($z_D^*$).
The results indicate that the greedy solution is within 1.7\% of the optimal solution values across all problem instances, 
the greedy algorithm identifies the optimal solution in $28$ of the $40$ problem instances. 
Gurobi takes considerably longer to load and solve the problem instances than the greedy algorithm for larger problem instances, which provides an incentive for its use.


\begin{table}\caption{Solution times for various instances of the defender's problem subject to a knapsack constraint and the relative accuracy of the greedy algorithm.}\label{tab:algm}\footnotesize
\centering
\begin{tabular}{| c | c | c | c | c | c | c | c | c |} \hline
\textbf{Layers} & \textbf{Nodes per} & \textbf{Edges} & \textbf{Controls} & \textbf{Budget} & \textbf{OptDef} &  \textbf{Greedy} &\textbf{Gap} \\
$\ell$ &\textbf{Layer}, \textbf{$|N_\ell|$} & $|E|$ & $|M|$ & $\beta$ & \textbf{Time (s)} & \textbf{Time (s)} & $z_{h} / z_D^*$ \\ \hline 
%
5 & 15 & 930 & 10 & 5 & 3.07 & 0.078 & 0.983 \\ \hline
        5 & 15 & 930 & 20 & 10 & 2.62 & 0.236 & 0.998 \\ \hline
        10 & 10 & 920 & 10 & 5 & 2.47 & 0.237 & 0.999 \\ \hline
        10 & 10 & 920 & 20 & 10 & 2.26 & 0.277 & 0.999 \\ \hline
        10 & 15 & 2055 & 12 & 6 & 19.9 & 0.998 & 0.997 \\ \hline
        10 & 15 & 2055 & 24 & 12 & 6.99 & 1.07 & 0.999 \\ \hline
        10 & 20 & 3640 & 14 & 7 & 49.7 & 2.25 & 0.998 \\ \hline
        10 & 20 & 3640 & 30 & 15 & 19.1 & 3.24 & 1.000 \\ \hline
        10 & 25 & 5675 & 16 & 8 & 66.2 & 2.12 & 0.999 \\ \hline
        10 & 25 & 5675 & 34 & 17 & 63.7 & 6.13 & 1.000 \\ \hline
        15 & 5 & 360 & 10 & 5 & 0.724 & 0.075 & 0.999 \\ \hline
        15 & 5 & 360 & 20 & 10 & 1.16 & 0.337 & 1.000 \\ \hline
        15 & 10 & 1420 & 12 & 6 & 8.40 & 0.328 & 0.997 \\ \hline
        15 & 10 & 1420 & 24 & 12 & 7.16 & 1.13 & 1.000 \\ \hline
        15 & 15 & 3180 & 14 & 7 & 24.0 & 1.13 & 0.993 \\ \hline
        15 & 15 & 3180 & 30 & 15 & 22.2 & 2.48 & 1.000 \\ \hline
        15 & 20 & 3000 & 16 & 8 & 22.1 & 1.15 & 1.000 \\ \hline
        15 & 20 & 3000 & 34 & 17 & 18.8 & 3.75 & 1.000 \\ \hline
        15 & 20 & 5640 & 16 & 8 & 68.1 & 2.36 & 1.000 \\ \hline
        15 & 20 & 5640 & 34 & 17 & 60.3 & 8.27 & 1.000 \\ \hline
        15 & 25 & 4000 & 18 & 9 & 36.7 & 2.22 & 0.999 \\ \hline
        15 & 25 & 4000 & 40 & 20 & 41.6 & 9.74 & 1.000 \\ \hline
        15 & 25 & 8800 & 18 & 9 & 157.6 & 4.15 & 1.000 \\ \hline
        15 & 25 & 8800 & 40 & 20 & 159.1 & 17.8 & 1.000 \\ \hline
        20 & 10 & 1920 & 14 & 7 & 12.0 & 0.717 & 1.000 \\ \hline
        20 & 10 & 1920 & 30 & 15 & 10.3 & 2.81 & 1.000 \\ \hline
        20 & 15 & 4000 & 16 & 8 & 38.7 & 1.33 & 1.000 \\ \hline
        20 & 15 & 4000 & 34 & 17 & 34.2 & 6.66 & 1.000 \\ \hline
        20 & 15 & 4305 & 16 & 8 & 44.8 & 1.62 & 1.000 \\ \hline
        20 & 15 & 4305 & 34 & 17 & 41.8 & 6.17 & 1.000 \\ \hline
        20 & 20 & 4000 & 18 & 9 & 12.1 & 0.637 & 1.000 \\ \hline
        20 & 20 & 4000 & 40 & 20 & 32.2 & 7.56 & 1.000 \\ \hline
        20 & 20 & 7640 & 18 & 9 & 124.6 & 3.54 & 1.000 \\ \hline
        20 & 20 & 7640 & 40 & 20 & 107.3 & 13.8 & 1.000 \\ \hline
        20 & 25 & 5000 & 20 & 10 & 22.7 & 1.09 & 1.000 \\ \hline
        20 & 25 & 5000 & 44 & 22 & 47.6 & 11.9 & 1.000 \\ \hline
        20 & 25 & 11925 & 20 & 10 & 95.1 & 2.83 & 1.000 \\ \hline
        20 & 25 & 11925 & 44 & 22 & 238.3 & 28.6 & 1.000 \\ \hline
        25 & 10 & 2420 & 16 & 8 & 5.96 & 0.459 & 1.000 \\ \hline
        25 & 10 & 2420 & 34 & 17 & 15.0 & 3.29 & 1.000 \\ \hline
\end{tabular}
\end{table}

\section{Conclusion}
\label{cyber:sec:conclusion}
In this paper, we introduce new models and algorithms for identifying a portfolio of security controls to deploy when considering multiple adaptive adversaries of varying levels of strategic sophistication.  To inform these investments, our approach extends an ARA framework to consider a maximum reliability path interdiction problem with a single defender and multiple attackers. All players are assumed to be boundedly rational, and we allow for uncertainty regarding system and player information. We present mixed integer programming formulations of the attacker and defender problems, and we introduce an ARA algorithm to iteratively solve the models. We also introduce an approximation algorithm that identifies near-optimal solutions to the defender's problem with a guaranteed $1-1/e$ approximation ratio. 

The solutions provide insight into investments that construct a layered security defense and that perform well against many adversaries, including some adversaries who are not strategic. One practical benefit of the proposed methodology is that it yields a suite of investment solutions, instead of a single solution, which can aid decision-makers. Another benefit of the modeling approach is that it allows for the consideration of non-strategic attackers, which could be used to model risks arises from nature, although we did not specifically consider the impact of natural disasters \citep{zhuang2007balancing}. \added{New vulnerabilities appear regularly, and therefore, organizations should proactively perform risk assessments to inform defensive investments and decision-making on a regular basis. The approach introduced in this paper seeks to help with these decisions.}

We illustrate the models and solution techniques on a case study. The framework indicates that the defensive strategies change along with the defender's level of strategic sophistication as well as those of the attackers. The results identify security controls that are effective across a range of adversarial assumptions. The solutions tend to defend better against less strategic attackers than more strategic defenders. Additionally, the results suggest that it may be better for the defender to overestimate rather than underestimate the strategic sophistication of the attackers.

There are several avenues for future research. First, the ARA approach with boundedly rational players can be extended to include other network interdiction models aside from the maximum reliability path interdiction problem. 
Second, model extensions could accommodate attackers with less knowledge of the network than the defender, possibly by eliminating portions of the network that the attacker does not know about when solving the attackers' formulation. 
Third, the ARA approach could be extended to balance the goal of security with system performance depending on the application under consideration.

\section*{Acknowledgments}
The work of the first and last authors was in part supported by the National Science Foundation Award 2000986. Any opinions, findings, and conclusions or recommendations expressed in this material are those of the authors and do not necessarily reflect the views of the National Science Foundation.
\added{The authors would also like to thank the anonymous reviewers, whose suggestions for improvement led to a substantially improved manuscript.}

\singlespacing
\bibliographystyle{plainnat}
\bibliography{Zotero}

\begin{thebibliography}{44}
\providecommand{\natexlab}[1]{#1}
\providecommand{\url}[1]{\texttt{#1}}
\expandafter\ifx\csname urlstyle\endcsname\relax
  \providecommand{\doi}[1]{doi: #1}\else
  \providecommand{\doi}{doi: \begingroup \urlstyle{rm}\Url}\fi

\bibitem[Ahuja et~al.(1993)Ahuja, Magnanti, and Orlin]{ahuja_network_1993}
Ravindra Ahuja, Thomas Magnanti, and James Orlin.
\newblock \emph{Network Flows: Theory, Algorithms, and Applications}.
\newblock Prentice-Hall, Inc., 1 edition, 1993.
\newblock ISBN 0-13-617549-X.

\bibitem[Albert et~al.(2023)Albert, Nikolaev, and Jacobson]{albert2023homeland}
Laura~A Albert, Alexander Nikolaev, and Sheldon~H Jacobson.
\newblock Homeland security research opportunities.
\newblock \emph{IISE Transactions}, 55\penalty0 (1):\penalty0 22--31, 2023.

\bibitem[Banks et~al.(2020)Banks, Gallego, Naveiro, and
  R{\'\i}os~Insua]{banks2020adversarial}
David Banks, V{\'\i}ctor Gallego, Roi Naveiro, and David R{\'\i}os~Insua.
\newblock Adversarial risk analysis: An overview.
\newblock \emph{Wiley Interdisciplinary Reviews: Computational Statistics},
  page e1530, 2020.

\bibitem[Bistarelli et~al.(2006)Bistarelli, Fioravanti, and
  Peretti]{bistarelli2006defense}
Stefano Bistarelli, Fabio Fioravanti, and Pamela Peretti.
\newblock Defense trees for economic evaluation of security investments.
\newblock In \emph{Proceedings of the First International Conference on
  Availability, Reliability and Security (ARES'06), Vienna, Austria, 20-22
  April 2006}, pages 8--pp. IEEE, 2006.

\bibitem[Caballero et~al.(2021)Caballero, Friend, and
  Blasch]{caballero2021adversarial}
William~N Caballero, Mark Friend, and Erik Blasch.
\newblock Adversarial machine learning and adversarial risk analysis in
  multi-source command and control.
\newblock \emph{Proceedings SPIE 11756, Signal Processing, Sensor/Information
  Fusion, and Target Recognition XXX}, 11756L:\penalty0 98--108, 2021.

\bibitem[Camerer et~al.(2004)Camerer, Ho, and Chong]{camerer_cognitive_2004}
C.~F. Camerer, T.-H. Ho, and J.-K. Chong.
\newblock A cognitive hierarchy model of games.
\newblock \emph{The Quarterly Journal of Economics}, 119\penalty0 (3):\penalty0
  861--898, 2004.

\bibitem[Cano et~al.(2016)Cano, Pollini, Falciani, and
  Turhan]{cano_modeling_2016}
Javier Cano, Alessandro Pollini, Lorenzo Falciani, and Uğur Turhan.
\newblock Modeling current and emerging threats in the airport domain through
  adversarial risk analysis.
\newblock \emph{Journal of Risk Research}, 19\penalty0 (7):\penalty0 894--912,
  2016.

\bibitem[Cavusoglu et~al.(2008)Cavusoglu, Raghunathan, and
  Yue]{cavusoglu_decision-theoretic_2008}
Huseyin Cavusoglu, Srinivasan Raghunathan, and Wei~T. Yue.
\newblock Decision-{{Theoretic}} and {{Game}}-{{Theoretic Approaches}} to {{IT
  Security Investment}}.
\newblock \emph{Journal of Management Information Systems}, 25\penalty0
  (2):\penalty0 281--304, 2008.

\bibitem[{Council of Economic
  Advisors}(2018)]{council_of_economic_advisors_cost_2018}
{Council of Economic Advisors}.
\newblock The cost of malicious cyber activity to the {U}.{S}. economy.
\newblock Technical report, The White House, Washington, D.C., 2018.
\newblock URL
  \url{https://www.whitehouse.gov/wp-content/uploads/2018/03/The-Cost-of-Malicious-Cyber-Activity-to-the-U.S.-Economy.pdf}.

\bibitem[DuBois(2020)]{dubois2020optimizations}
Eric DuBois.
\newblock \emph{Optimizations Models with Applications to Homeland Security
  Systems}.
\newblock PhD thesis, The University of Wisconsin-Madison, Madison, WI, 2020.

\bibitem[Fielder et~al.(2016)Fielder, Panaousis, Malacaria, Hankin, and
  Smeraldi]{fielder_decision_2016}
Andrew Fielder, Emmanouil Panaousis, Pasquale Malacaria, Chris Hankin, and
  Fabrizio Smeraldi.
\newblock Decision {{Support Approaches}} for {{Cyber Security Investment}}.
\newblock \emph{Decision Support Systems}, 86:\penalty0 13--23, 2016.

\bibitem[Hubbard and Seiersen(2016)]{hubbard2016measure}
Douglas~W Hubbard and Richard Seiersen.
\newblock \emph{How to measure anything in cybersecurity risk}.
\newblock John Wiley \& Sons, 2016.

\bibitem[Israeli and Wood(2002)]{israeli_shortest-path_2002}
Eitan Israeli and R.~Kevin Wood.
\newblock Shortest-path network interdiction.
\newblock \emph{Networks}, 40\penalty0 (2):\penalty0 97--111, 2002.

\bibitem[Joshi et~al.(2020)Joshi, Aliaga, and Insua]{joshi2020insider}
Chaitanya Joshi, Jesus~Rios Aliaga, and David~Rios Insua.
\newblock Insider threat modeling: An adversarial risk analysis approach.
\newblock \emph{IEEE Transactions on Information Forensics and Security},
  16:\penalty0 1131--1142, 2020.

\bibitem[Khuller et~al.(1999)Khuller, Moss, and Naor]{khuller1999budgeted}
Samir Khuller, Anna Moss, and Joseph~Seffi Naor.
\newblock The budgeted maximum coverage problem.
\newblock \emph{Information processing letters}, 70\penalty0 (1):\penalty0
  39--45, 1999.

\bibitem[Knowles et~al.(2015)Knowles, Prince, Hutchison, Disso, and
  Jones]{knowles_survey_2015}
William Knowles, Daniel Prince, David Hutchison, Jules Ferdinand~Pagna Disso,
  and Kevin Jones.
\newblock A survey of cyber security management in industrial control systems.
\newblock \emph{International Journal of Critical Infrastructure Protection},
  9:\penalty0 52--80, 2015.

\bibitem[Kruse et~al.(2017)Kruse, Frederick, Jacobson, and
  Monticone]{kruse_cybersecurity_2017}
Clemens~Scott Kruse, Benjamin Frederick, Taylor Jacobson, and D.~Kyle
  Monticone.
\newblock Cybersecurity in healthcare: {{A}} systematic review of modern
  threats and trends.
\newblock \emph{Technology and Health Care}, 25\penalty0 (1):\penalty0 1--10,
  2017.

\bibitem[Lallie et~al.(2020)Lallie, Debattista, and Bal]{lallie2020review}
Harjinder~Singh Lallie, Kurt Debattista, and Jay Bal.
\newblock A review of attack graph and attack tree visual syntax in cyber
  security.
\newblock \emph{Computer Science Review}, 35:\penalty0 100219, 2020.

\bibitem[Lee and Wolpert(2012)]{lee_game_2012}
Ritchie Lee and David~H. Wolpert.
\newblock Game theoretic modeling of pilot behavior during mid-air encounters.
\newblock In Tatiana~Valentine Guy, Miroslav Kárný, and David~H. Wolpert,
  editors, \emph{Decision {{Making}} with {{Imperfect Decision Makers}}.},
  volume~28 of \emph{Intelligent {{Systems Reference Library}}}. {Springer,
  Berlin, Heidelberg}, 2012.

\bibitem[Letchford and Vorobeychik(2013)]{Letchford-Vorobeychik2013}
J.~Letchford and Y.~Vorobeychik.
\newblock Optimal interdiction of attack plans.
\newblock In \emph{Proceedings of the 12th International Conference on
  Autonomous Agents and Multiagent Systems}, Saint Paul, MN, 2013.

\bibitem[Morton et~al.(2007)Morton, Pan, and Saeger]{morton_models_2007}
David~P. Morton, Feng Pan, and Kevin~J. Saeger.
\newblock Models for nuclear smuggling interdiction.
\newblock \emph{{IIE} Transactions}, 39\penalty0 (1):\penalty0 3--14, 2007.

\bibitem[Nandi et~al.(2016)Nandi, Medal, and Vadlamani]{nandi2016interdicting}
Apurba~K Nandi, Hugh~R Medal, and Satish Vadlamani.
\newblock Interdicting attack graphs to protect organizations from cyber
  attacks: A bi-level defender--attacker model.
\newblock \emph{Computers \& Operations Research}, 75:\penalty0 118--131, 2016.

\bibitem[{National Institute of Standards and Technology}(2018)]{nist2018}
{National Institute of Standards and Technology}.
\newblock Framework for improving critical infrastructure cybersecurity.
\newblock {NIST} report, National Institute of Standards and Technology,
  Washington, D.C., 2018.

\bibitem[Nemhauser et~al.(1978)Nemhauser, Wolsey, and
  Fisher]{nemhauser_analysis_1978}
George~L Nemhauser, Laurence~A Wolsey, and Marshall~L Fisher.
\newblock An analysis of approximations for maximizing submodular set
  functions—{I}.
\newblock \emph{Mathematical programming}, 14\penalty0 (1):\penalty0 265--294,
  1978.

\bibitem[Rios and Rios~Insua(2012)]{rios2012adversarial}
Jesus Rios and David Rios~Insua.
\newblock Adversarial risk analysis for counterterrorism modeling.
\newblock \emph{Risk Analysis}, 32\penalty0 (5):\penalty0 894--915, 2012.

\bibitem[Rios~Insua et~al.(2021)Rios~Insua, Couce-Vieira, Rubio, Pieters,
  Labunets, and G.~Rasines]{insua_adversarial_2019}
David Rios~Insua, Aitor Couce-Vieira, Jose~A Rubio, Wolter Pieters, Katsiaryna
  Labunets, and Daniel G.~Rasines.
\newblock An adversarial risk analysis framework for cybersecurity.
\newblock \emph{Risk Analysis}, 41\penalty0 (1):\penalty0 16--36, 2021.

\bibitem[Ross et~al.(2021)Ross, Pillitteri, Dempsey, Riddle, and
  Guissanie]{NIST800}
Ron Ross, Victoria Pillitteri, Kelley Dempsey, Mark Riddle, and Gary Guissanie.
\newblock Protecting controlled unclassified information in nonfederal systems
  and organizations.
\newblock {NIST} special publication 800-171v2, {National Institute of
  Standards and Technology}, Washington, D.~C., 2021.

\bibitem[Rothschild et~al.(2012)Rothschild, McLay, and
  Guikema]{rothschild2012adversarial}
Casey Rothschild, Laura McLay, and Seth Guikema.
\newblock Adversarial risk analysis with incomplete information: A level-k
  approach.
\newblock \emph{Risk Analysis: An International Journal}, 32\penalty0
  (7):\penalty0 1219--1231, 2012.

\bibitem[Salmer{\'o}n(2012)]{salmeron_deception_2012}
Javier Salmer{\'o}n.
\newblock Deception tactics for network interdiction: A multiobjective
  approach.
\newblock \emph{Networks}, 60\penalty0 (1):\penalty0 45--58, 2012.

\bibitem[Scheibehenne et~al.(2010)Scheibehenne, Greifeneder, and
  Todd]{scheibehenne_can_2010}
Benjamin Scheibehenne, Rainer Greifeneder, and Peter~M Todd.
\newblock Can there ever be too many options? {A} meta-analytic review of
  choice overload.
\newblock \emph{Journal of Consumer Research}, 37\penalty0 (3):\penalty0
  409--425, 2010.

\bibitem[Schneier(1999)]{Schneier1999}
B.~Schneier.
\newblock Attack trees: Modeling security threats.
\newblock \emph{Dr. Dobb's Journal of Software Tools}, 24\penalty0
  (12):\penalty0 21--29, 1999.

\bibitem[Singhal and Ou(2017)]{singhal2017security}
A.~Singhal and X.~Ou.
\newblock Security risk analysis of enterprise networks using probabilistic
  attack graphs.
\newblock In L.~Wang, S.~Jajodia, A.~Singhal, A.~Singhal, and X.~Ou, editors,
  \emph{Network Security Metrics}. Springer, Heidelberg, Germany, 2017.

\bibitem[Smith and Song(2020)]{smith_survey_2019}
J~Cole Smith and Yongjia Song.
\newblock A survey of network interdiction models and algorithms.
\newblock \emph{European Journal of Operational Research}, 283\penalty0
  (3):\penalty0 797--811, 2020.

\bibitem[Stahl and Wilson(1995)]{stahl_players_1995}
Dale~O. Stahl and Paul~W. Wilson.
\newblock On players' models of other players: Theory and experimental
  evidence.
\newblock \emph{Games and Economic Behavior}, 10\penalty0 (1):\penalty0
  218--254, 1995.

\bibitem[Stevens et~al.(2020)Stevens, Dykstra, Everette, Chapman, Bladow,
  Farmer, Halliday, and Mazurek]{stevens2020compliance}
Rock Stevens, Josiah Dykstra, Wendy~Knox Everette, James Chapman, Garrett
  Bladow, Alexander Farmer, Kevin Halliday, and Michelle~L Mazurek.
\newblock Compliance cautions: Investigating security issues associated with us
  digital-security standards.
\newblock In \emph{Network and Distributed Systems Security (NDSS) Symposium},
  San Diego, CA, 2020.

\bibitem[Sviridenko(2004)]{sviridenko_note_2004}
Maxim Sviridenko.
\newblock A note on maximizing a submodular set function subject to a knapsack
  constraint.
\newblock \emph{Operations Research Letters}, 32\penalty0 (1):\penalty0 41--43,
  2004.

\bibitem[Wang and Banks(2011)]{wang2011network}
Shouqiang Wang and David Banks.
\newblock Network routing for insurgency: An adversarial risk analysis
  framework.
\newblock \emph{Naval Research Logistics}, 58\penalty0 (6):\penalty0 595--607,
  2011.

\bibitem[Wang et~al.(2019)Wang, Di~Maio, and Zio]{wang2019adversarial}
Wei Wang, Francesco Di~Maio, and Enrico Zio.
\newblock Adversarial risk analysis to allocate optimal defense resources for
  protecting cyber--physical systems from cyber attacks.
\newblock \emph{Risk Analysis}, 39\penalty0 (12):\penalty0 2766--2785, 2019.

\bibitem[Wang and Lu(2013)]{wang2013cyber}
Wenye Wang and Zhuo Lu.
\newblock Cyber security in the smart grid: Survey and challenges.
\newblock \emph{Computer networks}, 57\penalty0 (5):\penalty0 1344--1371, 2013.

\bibitem[Zhang et~al.(2018)Zhang, Zhuang, and Behlendorf]{zhang2018stochastic}
Jing Zhang, Jun Zhuang, and Brandon Behlendorf.
\newblock Stochastic shortest path network interdiction with a case study of
  {A}rizona-{M}exico border.
\newblock \emph{Reliability Engineering \& System Safety}, 179:\penalty0
  62--73, 2018.

\bibitem[Zheng and Albert(2019{\natexlab{a}})]{zheng2019interdiction}
Kaiyue Zheng and Laura~A Albert.
\newblock Interdiction models for delaying adversarial attacks against critical
  information technology infrastructure.
\newblock \emph{Naval Research Logistics}, 66\penalty0 (5):\penalty0 411--429,
  2019{\natexlab{a}}.

\bibitem[Zheng and Albert(2019{\natexlab{b}})]{zheng2019robust}
Kaiyue Zheng and Laura~A Albert.
\newblock A robust approach for mitigating risks in cyber supply chains.
\newblock \emph{Risk Analysis}, 39\penalty0 (9):\penalty0 2076--2092,
  2019{\natexlab{b}}.

\bibitem[Zheng et~al.(2019)Zheng, Albert, Luedtke, and
  Towle]{zheng_budgeted_2019}
Kaiyue Zheng, Laura~A Albert, James~R Luedtke, and Eli Towle.
\newblock A budgeted maximum multiple coverage model for cybersecurity planning
  and management.
\newblock \emph{IISE Transactions}, 51\penalty0 (12):\penalty0 1303--1317,
  2019.

\bibitem[Zhuang and Bier(2007)]{zhuang2007balancing}
Jun Zhuang and Vicki~M Bier.
\newblock Balancing terrorism and natural disasters—defensive strategy with
  endogenous attacker effort.
\newblock \emph{Operations Research}, 55\penalty0 (5):\penalty0 976--991, 2007.

\end{thebibliography}
\end{document}